\documentclass[12pt,reqno]{article}

\usepackage[usenames]{color}
\usepackage{amssymb}
\usepackage{amsmath}
\usepackage{amsthm}
\usepackage{amsfonts}
\usepackage{amscd}
\usepackage{graphicx}

\usepackage[colorlinks=true,
linkcolor=webgreen,
filecolor=webbrown,
citecolor=webgreen]{hyperref}

\definecolor{webgreen}{rgb}{0,.5,0}
\definecolor{webbrown}{rgb}{.6,0,0}

\usepackage{color}
\usepackage{fullpage}
\usepackage{float}

\usepackage{graphics}
\usepackage{latexsym}
\usepackage{epsf}
\usepackage{breakurl}

\newcommand{\seqnum}[1]{\href{https://oeis.org/#1}{\rm \underline{#1}}}

\begin{document}


\theoremstyle{plain}
\newtheorem{theorem}{Theorem}
\newtheorem{corollary}[theorem]{Corollary}
\newtheorem{lemma}[theorem]{Lemma}
\newtheorem{proposition}[theorem]{Proposition}

\theoremstyle{definition}
\newtheorem{definition}[theorem]{Definition}
\newtheorem{example}[theorem]{Example}
\newtheorem{conjecture}[theorem]{Conjecture}

\theoremstyle{remark}
\newtheorem{remark}[theorem]{Remark}

\begin{center}
\vskip 1cm{\LARGE\bf Further Extensions of Sury's Identity
}
\vskip 1cm
\large

Gregory Dresden\\
Washington \& Lee University \\
Lexington, VA 24450 \\
USA\\
\href{mailto:dresdeng@wlu.edu}{\tt dresdeng@wlu.edu} \\
\vskip .2 in
Xiaoya Gao\\ 
\"Ur\"umqi Bayi High School  \\
\"Ur\"umqi, Xinjiang \\
China\\
\href{mailto:lexie131@163.com}{\tt lexie131@163.com} \\

\end{center}

\vskip .2 in
\begin{abstract} 
The equation commonly known as Sury's identity is a deceptively simple summation formula that connects the Lucas numbers, Fibonacci numbers, and powers of two. Many authors have given extensions and generalizations over the years; in this paper, we take a different approach that  allows us to produce a good number of new summation formulas, all from elementary (but non-trivial) methods.
\end{abstract}

\section{Introduction}

The identity
\begin{equation}
2^{n+1}F_{n+1} = \sum_{i=0}^n 2^i L_i
            \label{e.Sury}
\end{equation}
which connects the Fibonacci numbers, the Lucas numbers, and the powers of 2, is often called Sury's identity thanks to his deceptively simple proof \cite{Sury} that takes the sum of two geometric series. Fifteen years before Sury's proof, Benjamin and Quinn \cite{BQ99} proved the same identity by counting the tilings of a bracelet with squares and dominoes. 

Of course, there are many variations on this identity. Two years after 
Sury's proof, Martinjak \cite{Mar} used a telescoping sum to prove that 
\begin{equation}\label{e.Martinjak}
(-1/2)^{n}F_{n+1} = \sum_{i=0}^n (-1/2)^i L_{i+1},
\end{equation}
and if we multiply both sides by $(-2)^n$ then this becomes the convolution identity
\begin{equation}\label{e1}
F_{n+1} = \sum_{i=0}^n (-2)^{n-i} L_{i+1}.
\end{equation}
Marques \cite{Marques} used a quick induction argument to  produce a slightly more complicated version of equation (\ref{e.Sury}) with $3^i$ instead of $2^i$, giving us 
\begin{equation}
3^{n+1}F_{n+1} = \sum_{i=0}^n 3^i \big(L_i + F_{i+1}\big),\label{e.Marques}
\end{equation}
and a clever rearrangement of sums by Edgar \cite{Edgar} produced the more general equation
\begin{equation}\label{e.Edgar}
t^{n+1} F_{n+1} = 
 \sum_{i=0}^n t^i \big(L_i + (t-2)F_{i+1}\big). 
\end{equation}
We note that Edgar's equation (\ref{e.Edgar}) will easily give us Sury's equation (\ref{e.Sury}) when $t=2$, and Marques' equation (\ref{e.Marques}) when $t=3$. It is a rather pleasant challenge to show that it also gives Martinjak's equation (\ref{e1}) when $t=-1/2$. For two different methods of proving equation (\ref{e.Edgar}), with colorful tilings and with generating functions, see Martinjack and Prodinger \cite{MP}.

Dafnis \cite{Dafnis} used a telescoping sum to give a variation on Edgar's equation (\ref{e.Edgar}), which in our notation would be 
\begin{equation}\label{e.Dafnis}
(-1)^{n} F_{n+1} = 
 \sum_{i=0}^n (-1)^i t^{n-i} \big(L_{i+1}  + (t-2)F_{i}\big). 
\end{equation}
Dafnis went on to prove a similar identity involving the tribonacci numbers.

Abd-Elhameed and Zeyada took Edgar's equation (\ref{e.Edgar}) even further by taking up generalized Fibonacci and Lucas sequences, defined as 
	\begin{equation} \label{AEZ}
		\begin{aligned}
U_{n+2}^{a,b} &= 
a U_{n+1}^{a,b}  + b U_{n}^{a,b} \qquad \mbox{with \ }
U_{0}^{a,b} =0, \ U_{1}^{a,b} =1, \\[1.2ex]
V_{n+2}^{a,b} &= 
a V_{n+1}^{a,b}  + b V_{n}^{a,b} \qquad \mbox{with \ }
V_{0}^{a,b} =2, \ V_{1}^{a,b} =a. 
		\end{aligned}
	\end{equation}
Their version of Sury's identity (for $a,b,t$ all non-zero) is 
	\begin{equation} \label{AEZ2}
t^{n+1}U_{n+1}^{a,b} = \frac{1}{a}
\sum_{i=0}^n 
t^i \big(V_i^{a,b} + (at-2)U_{i+1}^{a,b}\big),
     \end{equation}
and we can clearly see how this is a 
generalization of 
Edgar's equation (\ref{e.Edgar}).
Bhatnagar considered even more geneal sequences by replacing the constants $a,b$ in equation (\ref{AEZ}) with 
sequences $a_n, b_n$. 
For other approaches, see recent papers by 
Adegoke and Frontczak \cite{Adegoke},
Chung, Yao, and Zhou \cite{CYZ},
Dafnis, Philippou, and 
Livieris \cite{DPL},
De Micheli Vitturi \cite{DMV},
Kuhapatanakul and Thongsing \cite{KT},
Kwong \cite{Kwong},
Melham \cite{Melham},
 and 
Philippou and Dafnis \cite{PD}.

Despite all this, we have discovered a new way to generalize  Sury's identity. Every other approach has tried to replace the $2^i$ in Sury's equation (\ref{e.Sury})
with 
something like the $3^i$ in equation (\ref{e.Marques}) and then the  $t^i$ as seen in equation (\ref{e.Edgar}), which leads to an adjustment in the other term in the summation, from the simple $L_i$ in equation (\ref{e.Sury}) to the slightly more complex $L_i + F_{i+1}$ in equation (\ref{e.Marques}) to the rather complicated $L_i + (t-2)F_{i+1}$ in equation 
(\ref{e.Edgar}), culminating in the general expression 
$V_i^{a,b} + (at-2)U_{i+1}^{a,b}$ from Abd-Elhameed and Zeyada's equation (\ref{AEZ2}). 
By contrast, our approach allows us to select the sequence (the Fibonacci numbers, or the Pell numbers, or really any linear recurrence sequence of order 2) inside the summation, and also the initial position, and then  produce the particular ``$t$" to make the formula work out nicely.

This is best explained by a demonstration. Recall that  Sury and Martinjak began their sums in equations (\ref{e.Sury}) and (\ref{e.Martinjak}) with $L_0$ and $L_1$ respectively, but with our approach we can  begin with $L_2$ or $L_3$ or even $L_4$,  giving us the new formulas  
\begin{align*}
F_{n+1} = \sum_{i=0}^{n}(-2)^{n-i} L_{i+1}
&= \frac{1}{3}\sum_{i=0}^{n}(-1/3)^{n-i} L_{i+2} \\
&= \frac{1}{4}\sum_{i=0}^{n}(-3/4)^{n-i} L_{i+3} \\
&= \frac{1}{7}\sum_{i=0}^{n}(-4/7)^{n-i} L_{i+4}. 
\end{align*}
However, our personal favorites are these generalizations of Sury's formula (and again, these seem to be new),
\begin{equation} \label{e.favorites}
		\begin{aligned}
\sum_{i=0}^{n} (-2)^{i}L_{3i} &= 2 -(-2)^{n+1}F_{3n},\\
 4\sum_{i=0}^{n} 9^{i}L_{6i} &= 9^{n+1}F_{6n} + 8,\\
17\sum_{i=0}^{n} (-38)^{i}L_{9i} &= 34 -(-38)^{n+1}F_{9n}.
    \end{aligned}
\end{equation} 
These are  similar (in spirit but not in fact) to the formulas we can derive from Theorem 2.1 in Adegoke and Frontczak's recent paper \cite{Adegoke}, which give us (after some manipulation)
\begin{align*}
 \sum_{i=0}^{n} (1/2)^{i}L_{3i} &= F_{3(n+1)}/2^n,\\
 4\sum_{i=0}^{n} (1/9)^{i}L_{6i} &= F_{6(n+1)}/9^n,\\
 17\sum_{i=0}^{n} (1/38)^{i}L_{9i} &= F_{9(n+1)}/38^n.
\end{align*}

We can also produce this next collection of formulas with the Fibonacci and Lucas numbers (which, again, we believe are new): 
	\begin{equation} \label{e.steps}
		\begin{aligned}
F_{2n+2} &= 1\cdot \sum_{i=0}^{n} 1^{n-i} F_{2i+1} 
            = 1\cdot \sum_{i=0}^{n} (-1)^{n-i} L_{2i+1},\\
F_{3n+3} &= 2\cdot \sum_{i=0}^{n} 1^{n-i} F_{3i+1}
            = 2\cdot \sum_{i=0}^{n} (-3)^{n-i} L_{3i+1},\\
F_{4n+4} & = 3\cdot \sum_{i=0}^{n} 2^{n-i} F_{4i+1}
            = 3\cdot \sum_{i=0}^{n} (-4)^{n-i} L_{4i+1}.
\end{aligned}
\end{equation}

Moving on to the Pell numbers (which are found at \seqnum{A000129} in 
the On-Line Encyclopedia of
Integer Sequences (OEIS) \cite{oeis}), we have 
	\begin{equation} \label{e.Pellsteps}
		\begin{aligned}
P_{n+1} 
&= \frac{1}{2}\sum_{i=0}^{n}(-1/2)^{n-i} P_{i+2} \\
&= \frac{1}{5}\sum_{i=0}^{n}(-2/5)^{n-i} P_{i+3} \\
&= \frac{1}{12}\sum_{i=0}^{n}(-5/12)^{n-i} P_{i+4}. 
\end{aligned}
\end{equation}
We compare these to equation (23) in Abd-Elhameed and Zeyada's recent paper \cite{AEZ}, which gives us 
\begin{equation}\label{e.Pellsteps-1}
(1/3)^{n+1}Q_{n+1} = 1-\frac{2}{3}\sum_{i=0}^n
(1/3)^{i}P_{i-1},
\end{equation}
 where $Q_n$ represents the $n$th Pell-Lucas number 
 \seqnum{A001333} as defined later in equation (\ref{e.PQ}). 

Finally, if we think of the Pell numbers (with initial values $P_0 = 0$ and $P_1 = 1$ and  recurrence $P_n = 2P_{n-1} + P_{n-2}$) as the ``silver Fibonacci numbers", then the ``bronze Fibonacci numbers" would be the sequence \seqnum{A006190} with same initial values but with recurrence 
$B_n = 3B_{n-1} + B_{n-2}$. Here are some interesting (and mostly new) formulas for these bronze Fibonacci numbers, which remind us of the similar identities in equation (\ref{e.steps}):
\begin{equation} \label{e.Bronze.steps}
		\begin{aligned}
B_{2n+2} &= 3\cdot \sum_{i=0}^{n} 1^{n-i} B_{2i+1},\\
B_{3n+3} &= 10\cdot \sum_{i=0}^{n} 3^{n-i} B_{3i+1},\\
B_{4n+4} & = 33\cdot \sum_{i=0}^{n} 10^{n-i} B_{4i+1}.
\end{aligned}
\end{equation}

There are countless variations on these formulas. Again, what sets our approach apart from all the previous work is that our weighted sums can start at an arbitrary position in our sequence, and our weighted sums avoid the complicated summands as seen in  equations (\ref{e.Edgar}) and (\ref{AEZ2}). Our proofs  rely mostly on induction, along with generalized versions of Cassini's  and d'Ocagne's identities. 

\section{Main Results}
We begin with a  theorem associated with a very particular type of sequence.
\begin{theorem}\label{t1}
Let $(A_n)_{n \geq 0}$ be a second-order linear recurrence sequence with initial term $A_0 = 1$ and with recurrence
\[
    A_{n} = c_1A_{n-1} + c_2A_{n-2}.
\]
If we set $t=\big(c_1 - A_1\big)$ then so long as $t \not= 0$, we have 
\begin{equation}\label{e.t1}
A_{n+2} -  A_1 A_{n+1} = 
\Big( A_2 - A_1^2\Big)\sum_{i=0}^n t^{n-i} A_{i}. 
\end{equation}
\end{theorem}

\begin{proof}
We proceed by induction. 
For the base case $n=0$, we find that equation (\ref{e.t1}) simplifies to 
\[
A_2 -A_1A_1 = \Big(A_2 - A_1^2\Big) \cdot A_0,
\]
and since $A_0 = 1$ we are done. 

Next, we assume the identity holds for some fixed $n$, and we try to prove it for $n+1$. If we let $R$ represent the  right-hand side of equation (\ref{e.t1})  with $n+1$ in place of $n$, we have 
\[
R=(A_2 - A_1^2) \sum_{i=0}^{n+1} t^{n+1-i}A_i,
\]
and if we separate out that last term in the summation and then factor out $t$ then we have 
\[
R=(A_2 - A_1^2) \Big( A_{n+1} + t\sum_{i=0}^n t^{n-i}A_i\Big).
\]
We multiply through to obtain
\[
R=(A_2 - A_1^2) A_{n+1} + t\cdot (A_2 - A_1^2) \sum_{i=0}^n t^{n-i}A_i,
\]
and we use our induction hypothesis on that last term to give us
\[
R=  (A_2 - A_1^2)A_{n+1} + t\cdot (A_{n+2} - A_1A_{n+1}).
\]
Recalling that $t = c_1 - A_1$, the above expression becomes
\[
R=  (A_2 - A_1^2)A_{n+1} + (c_1 - A_1)\cdot (A_{n+2} - A_1A_{n+1}).
\]
We see that the $A_1^2A_{n+1}$ terms cancel, and so we now have 
\begin{equation}\label{t1.temp1}
  R= A_2A_{n+1} + (c_1 - A_1)\cdot A_{n+2} - c_1A_1A_{n+1}.
\end{equation}
When we replace that leading $A_2$ term with
\[
A_2 = c_1A_1 + c_2A_0 =  c_1A_1 + c_2 1 = c_1 A_1 + c_2,
\]
then equation (\ref{t1.temp1}) becomes
\[
R= 
(c_1A_1 + c_2)A_{n+1} + (c_1 - A_1)\cdot A_{n+2} - c_1A_1A_{n+1}.
\]
Further simplification reduces this to 
\[
R= 
c_2A_{n+1} + (c_1 - A_1)\cdot A_{n+2}.
\]
Since $c_1 A_{n+2} + c_2A_{n+1}$ is equal to $A_{n+3}$, then we have 
\[
R= 
A_{n+3} - A_1 A_{n+2},
\]
and this is exactly equal to the left-hand side of equation (\ref{e.t1}) with $n+1$ in place of $n$, thus completing our induction step. 
\end{proof}

With all that taken care of, we are now ready for  our next theorem, which covers a much more general collection of sequences. 
\begin{theorem}\label{t2}
For $X_n$ a second-order recurrence sequence with initial values $X_0$ and $X_1$ and with recurrence 
\begin{equation}\label{e.sors}
X_n = c_1X_{n-1} + c_2X_{n-2},
\end{equation}
then for any fixed  integer $k$ such that $X_k$ and $X_{k-1}$ are both non-zero, we have 
\begin{equation}\label{e.t2}
X_0X_{n+2} - X_1X_{n+1} = \frac{X_0X_2-X_1^2}{X_k} 
\sum_{i=0}^n
\left(-c_2X_{k-1}/X_k\right)^{n-i} {X_{i+k}}.
\end{equation}
\end{theorem}

We note that we allow $k$ to be negative. 
Given that we have fixed initial values $X_0$ and $X_1$ for our sequence, then we can always define $X_{-1}$ and $X_{-2}$ and so on by the identity 
\[
X_{n-2} = \frac{1}{c_2}\left(X_n - c_1X_{n-1}\right),
\]
which follows from the definition of our sequence in equation (\ref{e.sors}). Our requirement of a  {\em second-order} recurrence sequence tells us that 
$c_2 \not=0$, and thus $X_{-1}$ and $X_{-2}$ (and so on) are well-defined.

\begin{proof}[Proof of Theorem \ref{t2}]
For $k$ fixed, we begin by setting 
\begin{equation}\label{e.An}    
A_n = \frac{X_{n+k}}{X_k}.
\end{equation}
This gives us the new sequence $(A_n)_{n \geq 0}$ with initial term $A_0 = 1$ but with the same recurrence as the sequence $(X_n)_{n \geq 0}$ in the statement of the theorem.
Next, 
for convenience we will set 
\[
t = c_1 - A_1,
\]
which means that 
\[
t = c_1 - X_{k+1}/X_k = (c_1 X_k - X_{k+1})/X_k = -c_2X_{k-1}/X_k,
\]
and since $X_{k-1}$ is also assumed to be non-zero then $t \not= 0$. This means we can now apply Theorem \ref{t1} to the sequence $(A_n)_{n \geq 0}$,
giving us
\begin{equation}\label{e.t1.again}
A_{n+2} -  A_1 A_{n+1} = 
\Big( A_2 - A_1^2\Big)\sum_{i=0}^n t^{n-i} A_{i}, 
\end{equation}
and it remains to show that we can transform equation (\ref{e.t1.again}) into equation (\ref{e.t2}) in the statement of the theorem.

The first step, of course, is to use equation (\ref{e.An}) to re-write equation (\ref{e.t1.again}) as
\begin{equation}\label{e.t1.again2}
\frac{X_{n+k+2}}{X_k} -  \frac{X_{k+1}}{X_k} \frac{X_{n+k+1}}{X_k} = 
\Big( \frac{X_{k+2}}{X_k} -  \frac{X_{k+1}^2}{X_k^2}\Big)\sum_{i=0}^n t^{n-i} \frac{X_{k+i}}{X_k}, 
\end{equation}
and if we multiply through by $X_k^2$ we arrive at 
\begin{equation}\label{e.t1.again3}
X_{n+k+2}X_k - X_{k+1} X_{n+k+1} = 
\frac{ X_{k+2} X_k  -  X_{k+1}^2}{X_k}\sum_{i=0}^n t^{n-i} X_{k+i}.
\end{equation}
This is almost, but not quite, a match for equation (\ref{e.t2}). On the left of 
equation (\ref{e.t1.again3}) we will show that 
\begin{equation}\label{e.left}
X_{n+k+2}X_k - X_{k+1} X_{n+k+1} = (-c_2)^k (X_{n+2}X_0 - X_{n+1}X_1),
\end{equation}
and on the right of equation equation (\ref{e.t1.again3}) we will show that 
\begin{equation}\label{e.right}
X_{k+2} X_k  -  X_{k+1}^2 = (-c_2)^k (X_{2}X_0 - X_{1}^2).
\end{equation}
After canceling  $(-c_2)^k$ from both sides, we will obtain a perfect match for equation (\ref{e.t2}).

So, it remains to establish the identities 
(\ref{e.left}) and (\ref{e.right}) for all integer values of $k$. 
Equation (\ref{e.left})  is just a generalization of d'Ocagne's identity, and it follows immediately from taking the determinant of both sides of the matrix equation
\[
\begin{pmatrix}
X_{n+k+2} & X_{k+1}\\
X_{n+k+1} & X_{k} 
\end{pmatrix}
= 
\begin{pmatrix}
c_1 & c_2 \\
1 & 0
\end{pmatrix}^k \cdot
\begin{pmatrix}
X_{n+2} & X_{1}\\
X_{n+1} & X_{0} 
\end{pmatrix},
\]
which holds for $k$ positive, negative, or zero. 

As for 
equation (\ref{e.right}), this 
is a generalization of Cassini's identity, and we will prove it by induction (twice), once for $k>0$ and again for $k<0$. 
First, take $k$ to be positive. We note that 
equation (\ref{e.right}) is certainly true for $k=0$, 
and if we assume it to be true for some fixed 
 $k$, then at $k+1$ we can write
\begin{equation*}
X_{k+3} X_{k+1}  -  X_{k+2}^2 = \left( c_1 X_{k+2} + c_2 X_{k+1}\right) X_{k+1} - X_{k+2}^2,
\end{equation*}
and we re-arrange 
the right-hand side to get 
\begin{equation*}
X_{k+3} X_{k+1}  -  X_{k+2}^2 = c_2 X_{k+1}^2 + X_{k+2}\left(c_1 X_{k+1} - X_{k+2}\right). 
\end{equation*}
Since $X_{k+2} = c_1X_{k+1} + c_2X_{k}$,then the above equation becomes
\begin{equation*}
X_{k+3} X_{k+1}  -  X_{k+2}^2 = c_2 X_{k+1}^2 -c_2  X_{k+2} X_k = (-c_2) \left(X_{k+2}X_k - X_{k+1}^2\right),
\end{equation*}
and now our induction hypothesis comes into play for that last expression on the right to give us 
\begin{equation*}
X_{k+3} X_{k+1}  -  X_{k+2}^2 =  (-c_2) (-c_2)^k (X_{2}X_0 - X_{1}^2) 
= (-c_2)^{k+1} (X_{2}X_0 - X_{1}^2),
\end{equation*}
thus establishing equation (\ref{e.right}) for $k$ positive. 

Finally, take $k$ to be negative. 
We again note that 
equation (\ref{e.right}) is certainly true for $k=0$, 
and if we assume it to be true for some fixed 
 $k$, then at $k-1$ we can write
\begin{equation*}
X_{k+1} X_{k-1}  -  X_{k}^2 = X_{k+1}\left( \frac{1}{c_2} X_{k+1}  -\frac{c_1}{c_2} X_{k}\right) - X_{k}^2,
\end{equation*}
and we re-arrange 
the right-hand side to get 
\begin{equation*}
X_{k+1} X_{k-1}  -  X_{k}^2 = \frac{1}{c_2} X_{k+1}^2 - \frac{X_{k}}{c_2}\Big(c_1 X_{k+1} + c_2 X_{k}\Big). 
\end{equation*}
Since $X_{k+2} = c_1X_{k+1} + c_2X_{k}$,then the above equation becomes
\begin{equation*}
X_{k+1} X_{k-1}  -  X_{k}^2 = \frac{1}{c_2} X_{k+1}^2 -\frac{X_k}{c_2} \Big( X_{k+2} \Big) = (-1/c_2) \left(X_{k+2}X_k - X_{k+1}^2\right),
\end{equation*}
and now our induction hypothesis comes into play for that last expression on the right to give us 
\begin{equation*}
X_{k+1} X_{k-1}  -  X_{k}^2 =  (-1/c_2) (-c_2)^k (X_{2}X_0 - X_{1}^2) 
= (-c_2)^{k-1} (X_{2}X_0 - X_{1}^2),
\end{equation*}
thus establishing equation (\ref{e.right}) for $k$ negative.
\end{proof}

\section{Examples}

We now show how to use Theorems \ref{t1} and  \ref{t2} to produce many of the identities from the Introduction. 

Here is an overview of our approach: we begin with a particular sequence (say, the Lucas numbers or the Pell numbers), and then we try out different variations and offsets so that we can apply one of our two theorems. Once we have obtained a generic summation identity from our theorem, we can then attempt various simplifications in order to arrive at an aesthetically-pleasing formula. For example, we might start with the Lucas numbers $L_n$, but then decide to apply our theorems with either $L_{n-1}$ (giving us Sury's identity in section \ref{s.31} below) or  $L_n$ (as seen in section 
\ref{s.32}, with various offsets) 
or $L_{jn}$ (as done in Theorem \ref{t.t4}) or even $L_{jn+k}/L_{k}$ (which we use 
to give us equation (\ref{e.gen.offsetL})
in Theorem \ref{t.t5}).

\subsection{Sury's identity}\label{s.31}
To prove that 
\[
2^{n+1}F_{n+1} = \sum_{i=0}^n 2^{i} L_{i},
\]
we will use Theorem \ref{t2} with $X_n = L_{n-1}$, which gives us  $X_0 = -1$, $X_1 = 2$, and $c_1 = c_2 = 1$. Taking $k=1$ in Theorem \ref{t2}, we get 
\begin{equation}
-L_{n+1} - 2L_{n} = \frac{-1-2^2}{2} 
\sum_{i=0}^n
\left(1/2\right)^{n-i} {L_{i}}.
\end{equation}
After cancelling the negatives, and re-writing $L_{n+1} + 2L_n = L_{n+2} + L_n = 5F_{n+1}$, we have 
\begin{equation}
5F_{n+1} = \frac{5}{2} 
\sum_{i=0}^n
\left(1/2\right)^{n-i} {L_{i}}.
\end{equation}
Finally, we cancel the $5$'s and we multiply both sides by $2^{n+1}$, to give us our desired formula, 
\[
2^{n+1}F_{n+1} = \sum_{i=0}^n 2^{i} L_{i}.
\]

\subsection{Martinjak's convolution and its generalizations}\label{s.32}
From equation (\ref{e1}) in the introduction we have 
\begin{equation*}
F_{n+1} = \sum_{i=0}^n (-2)^{n-i} L_{i+1},
\end{equation*}
which is part of the family of identities 
\begin{align*}
F_{n+1} = \sum_{i=0}^{n}(-2)^{n-i} L_{i+1}
&= \frac{1}{3}\sum_{i=0}^{n}(-1/3)^{n-i} L_{i+2} \\
&= \frac{1}{4}\sum_{i=0}^{n}(-3/4)^{n-i} L_{i+3} \\
&= \frac{1}{7}\sum_{i=0}^{n}(-4/7)^{n-i} L_{i+4} 
\end{align*}
from the introduction. The general form is
\begin{equation}\label{19}  
F_{n+1} =  \frac{1}{L_k}\sum_{i=0}^{n}(-L_{k-1}/L_k)^{n-i} L_{i+k},
\end{equation}
and to produce all of these from Theorem \ref{t2} we take $X_n = L_n$ so that we have   $X_0 = 2$, $X_1 = 1$, and $c_1 = c_2 = 1$. As a result, equation (\ref{e.t2}) in Theorem \ref{t2} becomes
\begin{equation}
2L_{n+2} - L_{n+1} = \frac{2\cdot 3 -1^2}{L_k} 
\sum_{i=0}^n
\left(-L_{k-1}/L_k\right)^{n-i} {L_{i+k}}.
\end{equation}
Since 
$2L_{n+2} - L_{n+1} = 5F_{n+1}$, we easily produce equation (\ref{19}).
\ 

\subsection{Fibonacci and Lucas $j$-step numbers}
We define the  ``$j$-step numbers" of a sequence to be 
the new sequence made up of every $j$th term in the original sequence. For example, when we talk about the 
$3$-step Fibonacci numbers then we could be referring to 
$(F_{3n})_{n \geq 0}$ or perhaps 
$(F_{3n+1})_{n \geq 0}$ or even 
$(F_{3n+2})_{n \geq 0}$, depending on our choice of offset. We assume throughout this section that 
$j\geq 1$, and we note that if $j=1$ then we get the ``regular" Fibonacci (or Lucas) numbers without any skipping.  

\subsubsection{Fibonacci numbers with zero offset}\label{3.3.0}
For our first case, we will consider the $j$-step Fibonacci numbers such as 
$F_{3n}$ for $j=3$ or 
$F_{4n}$ for $j=4$, with no offset. Here is our result.
\begin{theorem} 
For $j\geq 1$ an integer, the $j$-step Fibonacci numbers satisfy
\begin{equation}\label{e.jstep.zero.offset}
L_j^{n+1}F_{j(n-1)} + (-1)^{j(n+1)}F_{2j} = \sum_{i=0}^{n} (-1)^{j(n+i)}L_j^{i}F_{ji}.
\end{equation}
\end{theorem}

\begin{proof}
We begin with the 
well-known formula \cite[p.~112]{Kfib}
\begin{equation}\label{e.Ruggles}
F_{a+b} = L_{b}F_a + (-1)^{b+1}F_{a-b},
\end{equation}
and if we replace $a$ with $j(n-1)$ and $b$ with $j$, we get 
\begin{equation}\label{e.c1c2Fib}
F_{jn} = L_jF_{j(n-1)} + (-1)^{j+1}F_{j(n-2)}.
\end{equation}

For $j\geq 1$ fixed, we will use Theorem \ref{t2}  with $X_n= F_{jn}$, and since Theorem \ref{t2} requires both $X_k$ and $X_{k-1}$ to be non-zero, we will set $k=2$. Thanks to equation (\ref{e.c1c2Fib}), we have 
$c_1 = L_j$ and $c_2 = (-1)^{j+1}$. Since $X_0 = 0$ and $X_1 = F_j$, then the expression on the left-hand side of equation (\ref{e.t2}) of Theorem \ref{t2} is 
\begin{equation}\label{e.25}
X_0 X_{n+2} - X_1X_{n+1} = 0-F_jF_{j(n+1)},
\end{equation}
and the expression in front of the sum on the right-hand side of equation (\ref{e.t2}) is 
\begin{equation}\label{e.26}
\frac{X_0 X_{2} - X_1^2}{X_2} = \frac{0-F_j^2}{\ \ F_{2j}} = \frac{-F_j}{\ \ L_{j}},
\end{equation}
where we have used the well-known identity $F_{2j} = F_jL_j$. Likewise, the expression inside the sum in equation (\ref{e.t2}) is 
\begin{equation}\label{e.27}
(-c_2X_{1}/X_2)^{n-i} X_{i+2} = ((-1)^j F_j/F_{2j})^{n-i} F_{j(i+2)} = ((-1)^j /L_j)^{n-i}F_{j(i+2)}.
\end{equation}
When we substitute equations (\ref{e.25}), 
(\ref{e.26}), and
(\ref{e.27}) into equation (\ref{e.t2}), and cancel the $-F_j$ term from each side, we get
\begin{equation*}
F_{j(n+1)} = \frac{1}{L_j} \sum_{i=0}^n ((-1)^j/L_j)^{n-i}F_{j(i+2)}.
\end{equation*}
This equation actually looks a bit nicer if we multiply both sides by $(L_j)^{n+1}$ and simplify, giving us 
\begin{equation}
L_j^{n+1}F_{j(n+1)} = \sum_{i=0}^n (-1)^{j(n+i)}L_j^{i}F_{j(i+2)}.
\end{equation}
Further simplification is needed. First,  we replace $i$ with $i-2$ in the right-hand summation, giving us 
\begin{equation}
L_j^{n+1}F_{j(n+1)} = \sum_{i=2}^{n+2} (-1)^{j(n+i)}L_j^{i-2}F_{ji}.
\end{equation}
Then, we multiply both sides by $L_j^2$ to give us 
\begin{equation}
L_j^{n+3}F_{j(n+1)} = \sum_{i=2}^{n+2} (-1)^{j(n+i)}L_j^{i}F_{ji}.
\end{equation}
Next, we add in the ``missing term" for $i=1$ to both sides, giving us 
\begin{equation}
L_j^{n+3}F_{j(n+1)} + (-1)^{j(n+1)}L_jF_j = \sum_{i=1}^{n+2} (-1)^{j(n+i)}L_j^{i}F_{ji}.
\end{equation}
Finally, we replace $n$ with $n-2$, we replace $L_jF_j$ with $F_{2j}$, and we start the summation at $i=0$ instead of at $i=1$ to give us
\begin{equation}
L_j^{n+1}F_{j(n-1)} + (-1)^{j(n+1)}F_{2j} = \sum_{i=0}^{n} (-1)^{j(n+i)}L_j^{i}F_{ji},
\end{equation}
which is equation (\ref{e.jstep.zero.offset}) as desired.
\end{proof}

Here are the results for small values of $j$:
\begin{align*}
\mbox{(at $j=2$):}& \qquad 3^{n+1} F_{2(n-1)} + 3 = \sum_{i=1}^n 3^{i} F_{2i},\\
\mbox{(at $j=3$):}& \qquad 4^{n+1} F_{3(n-1)} + 8(-1)^{n+1} = (-1)^n\sum_{i=1}^n (-4)^{i} F_{3i},\\
\mbox{(at $j=4$):}& \qquad 7^{n+1} F_{4(n-1)} + 21 = \sum_{i=1}^n 7^{i} F_{4i},
\end{align*}
and it is interesting to compare them to this similar collection of formulas from Dresden and Tulskikh \cite[p.~2]{DT},
\begin{align*}
3 -  F_{2(n+2)}/3^n  &= \sum_{i=0}^n (1/3)^i F_{2i},\\
F_{3(n+2)}/4^n -8 &= \sum_{i=0}^n (1/4)^i F_{3i},\\
21-F_{4(n+2)}/7^n &= \sum_{i=0}^n (1/7)^i F_{4i}.
\end{align*}

\subsubsection{Lucas numbers with zero offset}\label{3.3.0b}

For the $j$-step Lucas numbers such as $L_{2n}$ for $j=2$ or $L_{3n}$ for $j=3$, 
we will again use Theorem \ref{t2}  but this time with 
$X_n= L_{jn}$ and with $k=1$. Here is our result.
\begin{theorem} \label{t.t4}
For $j\geq 1$ an integer, the $j$-step Lucas numbers satisfy
\begin{equation}\label{e.Lzero}
2(L_j/2)^{n+1}(F_{jn}/F_j) + 2(-1)^{jn}= \sum_{i=0}^{n} (-1)^{j(n+i)}(L_j/2)^{i}L_{ji},
\end{equation}
\end{theorem}

\begin{proof}
To begin, we will need the formula
\begin{equation}\label{e.c1c2Luc.part1}
L_{a+b} = L_{b}L_a + (-1)^{b+1}L_{a-b},
\end{equation}
which follows from adding together two copies of equation 
(\ref{e.Ruggles}), the first with $a+1$ in place of $a$, and the second with $a-1$ in place of $a$. If we replace $a$ with $j(n-1)$ and $b$ with $j$ in equation  (\ref{e.c1c2Luc.part1}), we get 
\begin{equation}\label{e.c1c2Luc.part2}
L_{jn} = L_{j}L_{j(n-1)} + (-1)^{j+1}L_{j(n-2)}.
\end{equation}
We set $X_n = L_{jn}$, which means that  $X_0 = 2$ and  $X_1 = L_j$. Thanks to equation (\ref{e.c1c2Luc.part2}) we have $c_1 = L_j$ and $c_2 = (-1)^{j+1}$. 

We now turn our attention to   the expression on the left-hand side of equation (\ref{e.t2}) of Theorem \ref{t2}, which  is 
\begin{equation}\label{e.37}
X_0 X_{n+2} - X_1X_{n+1} = 2L_{j(n+2)} -L_jL_{j(n+1)}.
\end{equation}
From equation (\ref{e.c1c2Luc.part2}) with $n+2$ in place of $n$ we have
	\begin{align*}
L_{j(n+2)} &= L_{j}L_{j(n+1)} + (-1)^{j+1}L_{jn},\\
	   \intertext{and from Exercise 55 in \cite[p.~118]{K} we have }
L_{j(n+2)} &= 5F_jF_{j(n+1)} -(-1)^{j+1}L_{jn},	\\
	\intertext{and so if we add these together we will get} 
2L_{j(n+2)} &= L_{j}L_{j(n+1)} + 5F_jF_{j(n+1)}.
\end{align*}
When we substitute this into the right-hand side of equation (\ref{e.37}) we get 
\begin{equation}\label{e.37b}
X_0 X_{n+2} - X_1X_{n+1} = 2L_{j(n+2)} -L_jL_{j(n+1)}
= 5F_jF_{j(n+1)}.
\end{equation}

Likewise, the expression in front of the sum on the right-hand side of equation (\ref{e.t2}) is 
\begin{equation}\label{e.34}
\frac{X_0 X_{2} - X_1^2}{X_1} = \frac{2L_{2j} -L_j^2}{\ \ L_{j}} = \frac{5F_j^2}{L_{j}},
\end{equation}
which follows from equation (\ref{e.37b}) with $n=0$. 

Now, since we are taking $k=1$ in this section, then the expression inside the sum in equation (\ref{e.t2}) is 
\begin{equation}\label{e.35}
(-c_2X_{0}/X_1)^{n-i} X_{i+1} = 
  ((-1)^j\cdot 2 /L_{j})^{n-i} L_{j(i+1)} .
\end{equation}
When we substitute equations (\ref{e.37b}), 
(\ref{e.34}), and
(\ref{e.35}) into equation (\ref{e.t2}), and cancel the $5F_j$ term from each side, we get
\begin{equation*}
F_{j(n+1)} = \frac{F_j}{L_{j}} \sum_{i=0}^n (2(-1)^j/L_j)^{n-i}L_{j(i+1)}.
\end{equation*}

This equation actually looks a bit nicer if we multiply both sides by $(L_j/2)^{n+1}$ and simplify, giving us 
\begin{equation}
(L_j/2)^{n+1}F_{j(n+1)} = \frac{F_j}{2}\sum_{i=0}^n (-1)^{j(n+i)}(L_j/2)^{i}L_{j(i+1)}.
\end{equation}

We will need to apply some further simplifications. First, we multiply both sides by $2/F_j$ and then replace $i$ with $i-1$ in the right-hand summation, giving us 
\begin{equation}
2(L_j/2)^{n+1}(F_{j(n+1)}/F_j) = \sum_{i=1}^{n+1} (-1)^{j(n+i-1)}(L_j/2)^{i-1}L_{ji}.
\end{equation}
We multiply both sides by $(L_j/2)$ to produce
\begin{equation}
2(L_j/2)^{n+2}(F_{j(n+1)}/F_j) = \sum_{i=1}^{n+1} (-1)^{j(n+i-1)}(L_j/2)^{i}L_{ji}.
\end{equation}
Next, we add in the ``missing term" for $i=0$ to  both sides, giving us 
\begin{equation}
2(L_j/2)^{n+2}(F_{j(n+1)}/F_j) + 2(-1)^{j(n+1)}= \sum_{i=0}^{n+1} (-1)^{j(n+i-1)}(L_j/2)^{i}L_{ji}.
\end{equation}
Finally, we replace $n$ with $n-1$ to get 
\begin{equation}
2(L_j/2)^{n+1}(F_{jn}/F_j) + 2(-1)^{jn}= \sum_{i=0}^{n} (-1)^{j(n+i)}(L_j/2)^{i}L_{ji},
\end{equation}
and this is equation (\ref{e.Lzero}), as desired. 
\end{proof}

We note that equation  (\ref{e.Lzero}) looks rather pleasant when $L_j$ is even; this happens when $j$ is a multiple of $3$. For our first such case,  we substitute $j=3$ into equation (\ref{e.Lzero}) and multiply both sides by $(-1)^n$ and simplify to get
\begin{equation}\label{e.Lzero.3}
2 -(-2)^{n+1}F_{3n} = \sum_{i=0}^{n} (-2)^{i}L_{3i}.
\end{equation}
For the second, we use $j=6$ and multiply both sides by $4$ and simplify to get 
\begin{equation}\label{e.Lzero.6}
9^{n+1}F_{6n} + 8 = 4\sum_{i=0}^{n} 9^{i}L_{6i}.
\end{equation}
These are the first two identities from equation (\ref{e.favorites}) in the introduction, and they are  quite nice and quite unexpected. (A similar procedure for $j=9$ or $j=12$ will give us additional equations.)

\subsubsection{Non-zero offsets}

The formulas with offsets $-1$, $+1$, and $+2$, when written down together, seem remarkably similar:
\begin{align*}
\mbox{(offset $-1$):} \qquad F_{j(n+1)}  & = F_j  \sum_{i=0}^{n} (F_{j+1})^{n-i} F_{ji-1} \qquad \mbox{for $j>1$},\\
\mbox{(offset $+1$):} \qquad F_{j(n+1)}  &= F_j  \sum_{i=0}^{n} (F_{j-1})^{n-i} F_{ji+1} \qquad \mbox{for $j>1$},\\
\mbox{(offset $+2$):} \qquad F_{j(n+1)}  &= F_j  \sum_{i=0}^{n} (-F_{j-2})^{n-i} F_{ji+2} \qquad \mbox{for $j>2$}.
\end{align*}
Here are two additional offset formulas, and there does indeed seem to be a pattern:
\begin{align*}
\mbox{(offset $+3$):} \qquad F_{j(n+1)}  &= \frac{F_j}{2}  \sum_{i=0}^{n} (F_{j-3}/2)^{n-i} F_{ji+3} \qquad \mbox{for $j>3$}, \\
\mbox{(offset $+4$):} \qquad F_{j(n+1)}  &= \frac{F_j}{3}  \sum_{i=0}^{n} (-F_{j-4}/3)^{n-i} F_{ji+4} \qquad \mbox{for $j>4$}. \\
\end{align*}
Here is the general form:

\begin{theorem}\label{t.t5}
For $j>|k|$ integers, the $j$-step Fibonacci numbers with offset $k$ satisfy 
\begin{equation}\label{e.gen.offsetF}
F_{j(n+1)}  = \frac{F_j}{F_k}  \sum_{i=0}^{n} ((-1)^{k+1}F_{j-k}/F_k)^{n-i} F_{ji+k} \qquad \mbox{for $j>|k|>0$,}
\end{equation}
and for the Lucas numbers we have 
\begin{equation}\label{e.gen.offsetL}
F_{j(n+1)}  = \frac{F_j}{L_k}  \sum_{i=0}^{n} ((-1)^{k}L_{j-k}/L_k)^{n-i} L_{ji+k} \qquad \mbox{for $j>|k|$.}
\end{equation}
\end{theorem}

We note that if we set $k=0$ in equation (\ref{e.gen.offsetL}) then we obtain Theorem 2.1 from Adegoke and Frontczak's recent paper \cite{Adegoke}, and if we set $k=1$ in equations (\ref{e.gen.offsetF})  and (\ref{e.gen.offsetL}) then we obtain the collection of identities in equation (\ref{e.steps}).

\begin{proof}
To establish   equation (\ref{e.gen.offsetF}), we will use  Theorem \ref{t1} with $A_n = F_{jn+k}/F_k$. Thus, 
we have $A_0 = 1$,  $A_1 = F_{j+k}/F_k$, and since equation (\ref{e.Ruggles}) tells us that 
\[
F_{jn+k} = L_jF_{j(n-1)+k} + (-1)^{j+1}F_{j(n-2)+k}, 
\]
then we have 
\[
A_{n} = L_jA_{n-1} + (-1)^{j+1}A_{n-2}. 
\]
Thus, for the conditions in the statement of Theorem \ref{t1} we have $c_1 = L_j$ and so  $t=c_1 - A_1 = L_j - F_{j+k}/F_k$, which gives us 
\[
t = L_j- F_{j+k}/F_k = \frac{L_jF_k - F_{j+k}}{F_k}, \]
and thanks again to equation (\ref{e.Ruggles}) which tells us that $F_{k+j} = L_j F_k +(-1)^{j+1}F_{k-j}$, this gives us 
\[
t = (-1)^{j}F_{k-j}/F_k.
\]
Since $k-j<0$, then we use the well-known fact that $F_{-n} = (-1)^{n+1}F_n$, giving us 
\begin{equation}\label{e.ttt}
t = (-1)^{j}(-1)^{k-j+1}F_{j-k}/F_k = 
(-1)^{k+1}F_{j-k}/F_k.
\end{equation}
We are almost ready to use 
equation (\ref{e.t1})
from Theorem \ref{t1}, but we still need to simplify the two expressions 
\[
A_{n+2} - A_1A_{n+1} \qquad \mbox{and} \qquad A_2 - A_1^2 
\]
that we find in 
equation (\ref{e.t1}). From our definition of $A_n$ as $F_{jn+k}/F_k$, these two expressions are 
\begin{equation}\label{e.cat1}
\frac{F_{j(n+2) + k}F_k - F_{j(n+1) + k}F_{j+k}}{F_k^2} 
    \qquad \mbox{and} \qquad 
\frac{F_{2j + k}F_k - F_{j + k}^2}{F_k^2} 
\end{equation}
Here we will need to call upon a generalization of Catalan's Identity 
\cite[p.~108]{Kfib}
which tells us that 
\begin{equation}\label{e.cat3}
 F_{a+c}F_{b-c} - F_{a}F_{b} = (-1)^{b+c+1} F_{a+c-b}F_{c}.
\end{equation}
For the first expression in equation (\ref{e.cat1}), we use equation (\ref{e.cat3}) with $a=j(n+1)+k$, $b=j+k$ and $c=j$ to give us 
\begin{equation}\label{e.cat4}
A_{n+2} - A_1A_{n+1} = \frac{F_{j(n+2) + k}F_k - F_{j(n+1) + k}F_{j+k}}{F_k^2}  =  \frac{(-1)^{k+1}F_{j(n+1)}F_j}{F_k^2},
\end{equation}
For the second expression in equation (\ref{e.cat1}), we simply take $n=0$ in equation (\ref{e.cat4}) to give us 
\begin{equation}\label{e.cat5}
A_2 - A_1^2 = \frac{F_{2j + k}F_k - F_{j + k}^2}{F_k^2} =  \frac{(-1)^{k+1}F_j^2}{F_k^2}.
\end{equation}

We can now use  Theorem \ref{t1} with $A_{n} = F_{jn+k}/{F_k}$, and using  our expressions in equations (\ref{e.ttt}), (\ref{e.cat4}), and (\ref{e.cat5}) in equation (\ref{e.t1}) we get 
\[
\frac{(-1)^{k+1}F_{j(n+1)}F_j}{F_k^2}
     =
\frac{(-1)^{k+1}F_j^2}{F_k^2}\sum_{i=0}^n 
\left( (-1)^{k+1}F_{j-k}/F_k \right)^{n-i}\frac{F_{ji+k}}{F_k}.
\]
With a bit of cancellation, this simplifies nicely  to give us 
\[
F_{j(n+1)}
     =
\frac{F_j}{F_k}\sum_{i=0}^n 
\left( (-1)^{k+1}F_{j-k}/F_k \right)^{n-i}F_{ji+k},
\]
which is  equation  (\ref{e.gen.offsetF}) as desired. 

To establish equation  (\ref{e.gen.offsetL}), we follow most of the same steps, but there are a few interesting new wrinkles. Proceeding as before, we will use Theorem \ref{t1} with $A_n = L_{jn+k}/L_k$. Thus, 
we have $A_0 = 1$,  $A_1 = L_{j+k}/L_k$, 
and since equation (\ref{e.c1c2Luc.part2}) tells us that 
\[
 L_{jn} = L_jL_{j(n-1)} + (-1)^{j+1}L_{j(n-2)},
\]
then we have 
\[
A_{n} = L_jA_{n-1} + (-1)^{j+1}A_{n-2}. 
\]
Thus, for the conditions in the statement of Theorem \ref{t1} we have $c_1 = L_j$ and so  $t=c_1 - A_1 = L_j - L_{j+k}/L_k$, which gives us 
\[
t = L_j- L_{j+k}/L_k = \frac{L_jL_k - L_{j+k}}{L_k}, \]
and thanks  to equation (\ref{e.c1c2Luc.part1}) which tells us that $L_{k+j} = L_j L_k +(-1)^{j+1}L_{k-j}$, this gives us 
\[
t = (-1)^{j}L_{k-j}/L_k.
\]
Since $k-j<0$, then we use the well-known fact that $L_{-n} = (-1)^{n}L_n$, giving us 
\begin{equation}\label{e.ttt.Luc}
t = (-1)^{j}(-1)^{k-j}L_{j-k}/L_k = 
(-1)^{k}L_{j-k}/L_k.
\end{equation}

We are almost ready to use 
equation (\ref{e.t1})
from Theorem \ref{t1}, but we still need to simplify the two expressions 
\[
A_{n+2} - A_1A_{n+1} \qquad \mbox{and} \qquad A_2 - A_1^2 
\]
that we find in 
equation (\ref{e.t1}). From our definition of $A_n$ as $L_{jn+k}/L_k$, these two expressions are 
\begin{equation}\label{e.cat1.Luc}
\frac{L_{j(n+2) + k}L_k - L_{j(n+1) + k}L_{j+k}}{L_k^2} 
    \qquad \mbox{and} \qquad 
\frac{L_{2j + k}L_k - L_{j + k}^2}{L_k^2} 
\end{equation}
Here we will need to produce a variant of Catalan's Identity 
but with Lucas numbers instead of Fibonacci numbers. 
If we take two versions of equation (\ref{e.cat3}) with $a+1$ and $a-1$ in place of $a$, we have
\begin{align*}
 F_{a+1+c}F_{b-c} - F_{a+1}F_{b} &= (-1)^{b+c+1} F_{a+1+c-b}F_{c}, \\
 F_{a-1+c}F_{b-c} - F_{a-1}F_{b} &= (-1)^{b+c+1} F_{a-1+c-b}F_{c},
\end{align*}
 and if we add these together we find that 
 \begin{equation}\label{e.halfFhalfL}
 L_{a+c}F_{b-c} - L_{a}F_{b} = (-1)^{b+c+1} L_{a+c-b}F_{c}.
 \end{equation}
We now repeat the process, but this time with equation 
(\ref{e.halfFhalfL}) and this time with $b+1$ and $b-1$ in place of $b$, giving us 
\begin{align*}
 L_{a+c}F_{b+1-c} - L_{a}F_{b+1} &= (-1)^{b+c} L_{a+c-b-1}F_{c}, \\
 L_{a+c}F_{b-1-c} - L_{a}F_{b-1} &= (-1)^{b+c} L_{a+c-b+1}F_{c}, 
\end{align*}
and if we add these together we find that 
 \begin{equation}\label{e.LL}
 L_{a+c}L_{b-c} - L_{a}L_{b} = 5(-1)^{b+c} F_{a+c-b}F_{c}.
 \end{equation}
This is not quite as beautiful as the Fibonacci version in equation (\ref{e.cat3}), but it will certainly meet our needs.

If we turn now to the first expression in equation (\ref{e.cat1.Luc}), we use equation (\ref{e.LL}) with $a=j(n+1)+k$, $b=j+k$ and $c=j$ to give us 
\begin{equation}\label{e.cat4.Luc}
A_{n+2} - A_1A_{n+1} = \frac{L_{j(n+2) + k}L_k - L_{j(n+1) + k}L_{j+k}}{L_k^2}  =  \frac{5(-1)^{k}F_{j(n+1)}F_j}{L_k^2},
\end{equation}
For the second expression in equation (\ref{e.cat1.Luc}), we simply take $n=0$ in equation (\ref{e.cat4.Luc}) to give us 
\begin{equation}\label{e.cat5.Luc}
A_2 - A_1^2 = \frac{L_{2j + k}L_k - L_{j + k}^2}{L_k^2} =  \frac{5(-1)^{k}F_j^2}{L_k^2}.
\end{equation}

We can now use  Theorem \ref{t1} with $A_{n} = L_{jn+k}/{L_k}$, and using  our expressions in equations (\ref{e.ttt.Luc}), (\ref{e.cat4.Luc}), and (\ref{e.cat5.Luc}) in equation (\ref{e.t1}) we get 
\[
\frac{5(-1)^{k}F_{j(n+1)}F_j}{L_k^2}
     =
\frac{5(-1)^{k}F_j^2}{L_k^2}\sum_{i=0}^n 
\left( (-1)^{k}L_{j-k}/L_k \right)^{n-i}\frac{L_{ji+k}}{L_k}.
\]
With a bit of cancellation, this simplifies nicely  to give us 
\[
F_{j(n+1)}
     =
\frac{F_j}{L_k}\sum_{i=0}^n 
\left( (-1)^{k}L_{j-k}/L_k \right)^{n-i}L_{ji+k},
\]
which is  equation  (\ref{e.gen.offsetL}) as desired. 
\end{proof}

\subsection{Pell and Pell-Lucas numbers} 

We recall \cite[p.~23]{K} that the Pell and Pell-Lucas numbers $P_n$ and $Q_n$ are defined as 
	\begin{equation} \label{e.PQ}
		\begin{aligned}
P_n  &= 
2P_{n-1} + P_{n-2} \qquad \mbox{with \ }
P_0 =0, \ P_1 =1, \\[1.2ex]
Q_n  &= 
2Q_{n-1} + Q_{n-2} \qquad \mbox{with \ }
Q_0 =1, \ Q_1 =1.	
        \end{aligned}
	\end{equation}
These are the sequences \seqnum{A000129} and \seqnum{A001333}, respectively. If we apply Theorem 2 with $X_n = P_n$ and with $k\not= 0,1$, then since $P_0P_{n+2} - P_1P_{n+1} = -P_{n+1}$
and 
$P_0P_2 - P_1^2 = -1$, then equation (\ref{e.t2}) gives us 
\[
P_{n+1} = \frac{1}{P_k} \sum_{i=0}^n (-P_{k-1}/P_k)^{n-i} P_{i+k}.
\]
If we take $k=2,3$, and $4$, we will produce the identities in equation (\ref{e.Pellsteps}). However, it is interesting to note that with $k=-1$ we produce the equation
\begin{equation}\label{e.PPtemp}
P_{n+1} =  \sum_{i=0}^n 2^{n-i} P_{i-1},
\end{equation}
which is quite different from equation (\ref{e.Pellsteps-1}). However, it is equivalent to the almost-identical formula in \cite[page 3]{DT},
\[
P_{n+2} = 2^{n+1} + \sum_{i=0}^n 2^{n-i}P_{i}.
\]
To see this, we multiply both sides of equation (\ref{e.PPtemp}) by 2,  we replace $i-1$ with $i$, we pull out the $n$th term, we replace $2P_{n+1}$ with $P_{n+2} - P_n$, and we absorb $P_n$ into the sum.

\subsection{The bronze Fibonacci numbers}

As mentioned in the introduction, we define the ``bronze Fibonacci numbers" as $B_0 = 0$, $B_1 = 1$, and 
\[
B_n = 3B_{n-1} + B_{n-2}.
\]
This gives us the sequence \seqnum{A006190}. To prove the first identity in equation (\ref{e.Bronze.steps}),
we will use Theorem \ref{t1} with $A_n = B_{2n+1}$, and so $A_0 = B_1 = 1$ and $A_1 = B_3 =10$. Starting with 
\[
B_n = 3B_{n-1} + B_{n-2},
\]
then it is easy to show that 
\[
B_{2n+1} = 11B_{2n-1} - B_{2n-3}
\]
which means that 
\[
A_{n} = 11A_{n-1} - A_{n-2}.
\]
So, using Theorem \ref{t1} with $t=c_1 - A_1 = 11-10 = 1$,  we have 
\begin{equation}\label{e.Bronze1}
A_{n+2} - 10 A_{n+1} = (A_2 - A_1^2)\sum_{i=0}^n 1^{n-1}A_i.
\end{equation}
Now, $A_{n+2}-10A_{n+1}$ simplifies as follows:
\begin{align*}
A_{n+2}-10 A_{n+1} &= (11A_{n+1} - A_n) -10 A_{n+1}\\
&= A_{n+1} - A_n \\
&= B_{2n+3} - B_{2n+1} \\
&= 3B_{2n+2}.
\end{align*}
Furthermore, $A_2 - A_1^2$ is equal to $B_5 - B_3^2$, and since $B_3 = 10$ and $B_5 = 109$, then putting everything together we see that equation (\ref{e.Bronze1}) becomes
\begin{equation}\label{e.Bronze2}
3B_{2n+2} = 9\sum_{i=0}^n 1^{n-1}B_{2i+1},
\end{equation}
and this simplifies nicely to give us the first identity in equation (\ref{e.Bronze.steps}). 
The other two identities have similar proofs. 
\section{Conclusion}

As we  mentioned in the introduction, there are countless other identities that we can prove: we simply select a sequence, choose a starting value, and then apply our theorem. To take a random example, if we let $a_n$ be the sequence with initial values $a_0 = 0$ and $a_1 = 1$ and with recurrence $a_n = 4a_{n-1} + 3a_{n-2}$, then we can use Theorem \ref{t2} to quickly produce the following identities for this sequence (which is sequence 
\seqnum{A015530} in the OEIS):
\[
a_{n+1} 
 = \frac{1}{4}\sum_{i=0}^n (-3/4)^{n-i}a_{i+2} 
 = \frac{1}{19}\sum_{i=0}^n (-12/19)^{n-i}a_{i+3}. 
\]

We invite the reader to try their hand at developing new such identities. It would also be reasonable (and quite interesting) to try adapting  the methods of this paper to the Fibonacci and Lucas {\em polynomials}, as seen in Adegoke and Frontczak's recent article \cite{Adegoke}. The possibilities are unlimited.

\section{Acknowledgments}
We are grateful to Pioneer Academics for arranging the summer research project which led to this paper.

\bigskip
\hrule
\bigskip

\noindent 2020 {\it Mathematics Subject Classification}: Primary 11B39; 
Secondary 11B37. 

\noindent \emph{Keywords: }  Sury's identity, Fibonacci numbers, Lucas numbers, convolution, weighted sum.

\bigskip
\hrule
\bigskip

\noindent (Concerned with sequences
\seqnum{A000032}, 
\seqnum{A000045}, 
\seqnum{A000129},
\seqnum{A001333} 
\seqnum{A006190},
and
\seqnum{A015530}.)

\end{document}